\documentclass[11pt,oneside,a4paper]{article}
\usepackage[latin1]{inputenc}
\usepackage{ulem}
\usepackage{color}
\usepackage{fancybox}
\usepackage{fancyhdr}
\usepackage{amssymb}
\usepackage{amsmath}
\usepackage{amsthm}
\usepackage{stmaryrd}
\usepackage{graphicx}
\usepackage{hyperref}
\usepackage{epic}
\usepackage{fixmath}
\usepackage{bbm}
\usepackage{supertabular}
\usepackage{longtable}
\usepackage[left=3cm,right=3cm,top=2.5cm,bottom=2cm]{geometry}

\newcommand{\F}{\mathbb{F}}
\newcommand{\HH}{\mathbb{H}}
\newcommand{\C}{\mathbb{C}}
\newcommand{\R}{\mathbb{R}}

\newcommand{\Z}{\mathbb{Z}}

\newcommand{\GL}{\mathrm{GL}}
\newcommand{\SL}{\mathrm{SL}}
\newcommand{\Sp}{\mathrm{Sp}}
\newcommand{\U}{\mathrm{U}}
\newcommand{\SU}{\mathrm{SU}}
\newcommand{\SO}{\mathrm{SO}}
\newcommand{\Orm}{\mathrm{O}}
\newcommand{\Spin}{\mathrm{Spin}}
\newcommand{\Span}{\mathrm{span}}
\newcommand{\Phirm}{\mathrm{\Phi}}
\newcommand{\Iso}{\mathrm{Iso}}
\newcommand{\Mink}{\mathrm{Mink}}
\newcommand{\g}{\mathfrak{g}}
\newcommand{\Ad}{\mathrm{Ad}}

\newtheorem{theorem}{Theorem}[section]
\newtheorem{proposition}{Proposition}[section]
\newtheorem{lemma}{Lemma}[section]
\newtheorem{corollary}{Corollary}[section]
\newtheorem{remark}{Remark}[section]

\font\tenrm=cmr10
\font\cmssl=cmss10 at 12 pt
\font\bigss=cmssdc10 scaled 2300

\begin{document}
\rightline{}
\vskip 1.5 true cm
\begin{center}
{\bigss Pseudo-Riemannian almost hypercomplex homogeneous spaces with irreducible isotropy}
\vskip 1.0 true cm
{\cmssl  V.\ Cort\'es$^1$ and B.\ Meinke$^2$
} \\[3pt]
{\tenrm $^1$Department of Mathematics\\
and Center for Mathematical Physics\\
University of Hamburg\\
Bundesstra{\ss}e 55,
D-20146 Hamburg, Germany\\
vicente.cortes@uni-hamburg.de 
\\[1em]
$^2$Department of Mathematics\\ 
University of D\"usseldorf\\
Universit\"atsstra{\ss}e 1, Raum 25.22.03.62 \\
D-40225 D\"usseldorf , Germany\\  
Benedict.Meinke@uni-duesseldorf.de}
\end{center}
\vskip 1.0 true cm
\baselineskip=18pt
\begin{abstract}
\noindent
We classify homogeneous pseudo-Riemannian manifolds of index $4$ which admit an invariant almost hyper-Hermitian structure and an $\HH$-irreducible isotropy group. The main result is that all these spaces are hyper-K\"ahler and flat except in dimension $12$.
\end{abstract}

\section{Introduction}
Ahmed and Zeghib \cite{AhmedZeghib} investigated homogeneous almost complex manifolds of index $2$ which they call Hermite-Lorentz spaces referring to the complex index $1$. More precisely they considered a manifold $M=G/H$ of dimension $2n+2\geq 8$ where $G\subset\mathrm{Iso}(M,g,J)$ is a connected Lie group that acts transitively on $M$ and $\mathrm{Iso}(M,g,J)$ denotes the subgroup of the isometry group where elements preserve the almost complex structure. They showed that $M$ is a K\"ahler manifold locally isometric to one of the following spaces
$$ \Mink_{n+1}(\C)=\C^{1,n}, \ \mathrm{dS}_{n+1}(\C)=\SU(1,n+1)/\U(1,n), \ \mathrm{AdS}_{n+1}(\C)=\SU(2,n)/\U(1,n),$$
$$ \C\mathrm{dS}_{n+1}=\SO^0(1,n+2)/\SO^0(1,n)\times\SO(2), \ \C\mathrm{AdS}=\SO^0(3,n)/\SO(2)\times\SO^0(1,n).$$
Their strategy was to consider $H^0$ as a connected $\C$-irreducible Lie subgroup of $\U(1,n)$ and to use that such groups contain the subgroup $\SO^0(1,n)$. Since the K\"ahler form $\omega$ and the Nijenhuis tensor $N$ are respectively $H$-invariant and $H$-equivariant they are preserved under the action of $\SO^0(1,n)\subset H$. From this they conclude that $\omega$ is closed and that $N$ vanishes. Hence, the homogeneous spaces are already K\"ahler manifolds. A detailed investigation of the possibilities for $H^0$ gives the above list of manifolds.\\

In this paper we follow the idea of Ahmed and Zeghib but consider instead almost hypercomplex pseudo-Riemannian manifolds of index $4$ which 
have an $\HH$-irreducible isotropy group. It turns out that except in dimension $12$ all the manifolds are flat. Our main result is the following theorem.
\begin{theorem}\label{th:MainTheorem}
Let $(M,g,(J_1,J_2,J_3))$ be a connected almost hypercomplex pseudo\--Rieman\-nian manifold of index $4$ and $\dim M=4n+4\geq 8$, such that there exists a connected Lie subgroup $G\subset \Iso(M,g,(J_1,J_2,J_3))$ acting transitively on $M$. If the isotropy group $H=G_p$, $p\in M$, acts $\HH$-irreducibly, then $(M,g,(J_1,J_2,J_3))$ is globally hyper-K\"ahler and locally isometric to 
$\Mink_{n+1}(\HH)$ or $\dim M=12$ and $H^0$ is either $H^0\cong\SO^0(1,2)$, $H^0\cong\SU(1,2)$ or trivial.
\end{theorem}
Here $\mathrm{Iso}(M,g,(J_1,J_2,J_3))$ denotes the subgroup of $\mathrm{Iso}(M,g)$ the elements of which preserve the three almost complex structures $J_1$, $J_2$, and $J_3$. The quaternionic Minkowski space is denoted by $\Mink_{n+1}(\HH)=\HH^{1,n}$.\\

\noindent\textbf{Acknowledgments.} This work was partly supported by the German Science Foundation (DFG) under the
Collaborative Research Center (SFB) 676 Particles, Strings and the Early Universe. We would like to thank Ines Kath for valuable remarks.
\section{Hyperbolic spaces and subgroups of Sp(1,\textit{n})}
\subsection{Facts about hyperbolic spaces}
In this Section we collect some known facts about the quaternionic hyperbolic spaces which will be needed for the proof of our main result. We will use the same notation as in \cite{ChenGreenberg}. In the following $\F$ denotes the real numbers $\R$, the complex numbers $\C$ or the quaternions $\HH$. 
Recall first the classification of the totally geodesic submanifolds of the hyperbolic space $H^n(\F)$ over $\F$.
\begin{proposition}[{\cite[Prop. 2.5.1]{ChenGreenberg}}]\label{prop:totallygeodesicsubmanifolds}
Any totally geodesic submanifold of $H^n(\F)$ is equivalent under $\U(1,n;\F)$ to one of the following:
	\begin{itemize}
		\item[(i)] $H^m(\F')$, where $\F'\subseteq \F$, $\F'=\R,\C$ or $\HH$, and $1\leq m \leq n$;
		\item[(ii)] $H^1(\mathbb{I}):=e_1\mathbb{I}\cap B^n(\HH)$, where $e_1=(1,0,\ldots,0)^T\in\HH^n$ and $\mathbb{I}=\Span_\R\left\{\mathbold{i},\mathbold{j},\mathbold{k}\right\}$. It occurs only if $\F=\HH$.
	\end{itemize}
	These are all inequivalent under $\U(1,n;\F)$.
\end{proposition}
\begin{proposition}[{\cite[Proposition 4.2.1]{ChenGreenberg}}]\label{prop:decompositionstabilizer}
Let $M$ be a totally geodesic submanifold in $H^n(\F)$ and let $I(M)$ be the stabilizer of $M$ in $\U(1,n;F)$. Let $K(M)$ be the subgroup of $I(M)$ which leaves $M$ pointwise fixed. Then there exists a Lie subgroup $U(M)\subset I(M)$ such that $I(M)=K(M)U(M)$ (almost semidirect product). The identity component $U^0(M)$ is a simple Lie group when $\dim M>1$, and $I^0(M)=K^0(M)U^0(M)$ is an almost direct product. $U^0(M)$ induces the connected isometry group of $M$.
\end{proposition}
\noindent The following table covers all possibilities of $I(M)$ for a totally geodesic submanifold $M\subset H^n(\HH)$.\\[3pt]
{\footnotesize
\begin{longtable}{|l| l|}
\caption{Decomposition of $I(M)$}\\ \hline
   & $I(M)=K(M)U(M)$
\endhead
\hline
\label{table:Liesubgroups}
$M=H^m(\HH)\subset H^n(\HH)$ & $K(M)=\left\{\pm\mathbbm{1}_{m+1}\right\}\times \Sp(n-m)$, $U(M)=\Sp(1,m)\times \left\{\mathbbm{1}_{n-m}\right\}$,\\
                    & $K^0(M)=\left\{\mathbbm{1}_{m+1}\right\}\times\Sp(n-m)$;\\
										\hline
$M=H^m(\C)\subset H^n(\HH)$  & $K(M)=\U(1)\cdot\mathbbm{1}_{m+1}\times\Sp(n-m)$, $U(M)=\SU(1,m)\cdot\left\{\pm 1,\pm \mathbold{j}\right\}\times\left\{\mathbbm{1}_{n-m}\right\},$\\
                    & $U^0(M)=\SU(1,m)\times\left\{\mathbbm{1}_{n-m}\right\}$;\\
\hline
$M=H^m(\R)\subset H^n(\HH)$ & $K(M)=\Sp(1)\cdot\mathbbm{1}_{m+1}\times\Sp(n-m)$, $U(M)=\Orm(1,m)\times\left\{\mathbbm{1}_{n-m}\right\}$,\\
                   & $U^0(M)=\SO^0(1,m)\times\left\{\mathbbm{1}_{n-m}\right\}$;\\
\hline
$M=H^1(\mathbbm{I})\subset H^n(\HH)$ & $K(M)=\left\{\pm\mathbbm{1}_2\right\}\times\Sp(n-1)$, $U(M)=U\times\left\{\mathbbm{1}_{n-1}\right\}$,\\
                     & $K^0(M)=\left\{\mathbbm{1}_2\right\}\times\Sp(n-1)$, $U^0(M)=U^0\times\left\{\mathbbm{1}_{n-1}\right\}$.\\
\hline
\end{longtable}
}
\noindent In the case $M=H^1(\mathbbm{I})$ the Lie group $U\subset \Sp(1,1)$ is given by
$$ U= \left\{A\in \Sp(1,1)\left| A=\begin{pmatrix} a & -b\\ \varepsilon b & \varepsilon a \end{pmatrix}, \varepsilon=\pm 1\right.\right\}. $$
One can show that the elements of $U^0$ are precisely the elements of $\Sp(1,1)$ which commute with $$\Phirm=\begin{pmatrix} 0 & -1\\ 1 & 0 \end{pmatrix}.$$
Notice that $U^0$ acts $\HH$-irreducibly on $\HH^2$. It is sufficient to check this for the Lie algebra $\mathfrak{u}$. The matrices
$$ x=\frac{\mathbold{i}}{2}\Phirm,\ y=\frac{\mathbold{j}}{2}\Phirm,\ z=\frac{\mathbold{k}}{2}\Phirm, $$
$$ u=\frac{\mathbold{i}}{2}\mathbbm{1}_2,\ v=\frac{\mathbold{j}}{2}\mathbbm{1}_2, w=\frac{\mathbold{k}}{2}\mathbbm{1}_2 $$
form a basis of $\mathfrak{u}\subset\mathfrak{sp}(1,n)$. We have the following eigenspace decomposition for $x$
$$ \HH^{1,1} = \begin{pmatrix} \mathbold{i}\\ -1 \end{pmatrix}\cdot\HH\oplus \begin{pmatrix} \mathbold{i}\\ 1 \end{pmatrix}\cdot\HH.$$
Assume there exists a $\mathfrak{u}$-invariant subpace $V$ of quaternionic dimension one. Then $V$ is one of the quaternionic eigenspaces of $x$. But these spaces are not preserved by $y$.
Hence, $\mathfrak{u}$ acts $\HH$-irreducibly.\\
Furthermore, one can show that $U^0$ is simply connected and that its Lie algebra $\mathfrak{u}\cong\mathfrak{so}(1,3)$. This implies $U^0\cong \Spin^0(1,3)$.\\
\noindent Recall that the elements of $\Sp(1,n)$ are classified according to their fixed points in $\overline{H^n(\HH)}$. An element $g\in\Sp(1,n)$ is called \textbf{elliptic} if it has one fixed point in $H^n(\HH)$. It is called \textbf{parabolic} if it has exactly one fixed point and this point lies on the boundary $\partial H^n(\HH)$. If $g$ has exactly two fixed points which lie on the boundary it is called \textbf{loxodromic}. Any element with three or more fixed points on the boundary has also a fixed point in $H^n(\HH)$. Hence, the above classification covers all possibilities.\\
\noindent If $g\in \Sp(1,n)$ and $v\in \HH^{1,n}$, $\lambda\in\HH$, such that $g(v)=v\lambda$, then $A(v\mu)=(v\mu)\mu^{-1}\lambda\mu$ for all $\mu\in\HH\setminus\left\{0\right\}$, i.e.\ $v\mu$ is an eigenvector of $g$ with eigenvalue $\mu^{-1}\lambda\mu$. Thus the eigenvalues of $g$ occur in similarity classes. These are called of negative or positive type if the corresponding eigenvector is timelike or spacelike, respectively.
\begin{proposition}[{\cite[Prop. 3.2.2]{ChenGreenberg}}]\label{prop:ellipticfixedpoints}
Let $g\in\U(1,n;\F)$ be elliptic, let $\Lambda_0$ be its negative class of eigenvalues, and let $\Lambda_1,\ldots,\Lambda_n$ be its positive classes. Let $F(g)$ denote the set of fixed points of $g$ in $H^n(\F)$.
\begin{itemize}
	\item[(i)] If $\Lambda_0\neq \Lambda_i$ for all $1\leq i\leq n$, then $F(g)$ contains only one point.
	\item[(ii)] Suppose that $\Lambda_0$ coincides with exactly $m$ of the classes $\Lambda_i$, $1\leq i\leq n$. Then $F(g)$ is a totally geodesic submanifold, which is equivalent to $H^m(\F)$ if $\Lambda_0\subset\R$, and to $H^m(\C)$ if $\Lambda_0\not\subset\R$.
\end{itemize}
\end{proposition}
\begin{remark} One has to pay attention to the notation. The authors of \cite{ChenGreenberg} denote by $\C$ a subfield of $\F$ which contains $\R$ and is isomorphic to the field of complex numbers. Hence, in Proposition $\ref{prop:ellipticfixedpoints}$, $\C$ could be for example $\Span_\R\left\{1,\mathbold{j}\right\}$.
\end{remark}
\begin{lemma}[{\cite[Lemma 3.2.2]{ChenGreenberg}}]\label{le:loxodromicconjugate}
  Let $g\in \U(1,n;\F)$ which fixes $\pm f_1=(\pm 1,0,\ldots,0)^T\in \partial H^n(\F)$ (considered as elements of the sphere). Then $$ g=\begin{pmatrix} c\lambda & s\lambda & 0\\ s\lambda & c\lambda & 0\\ 0 & 0 & A \end{pmatrix}, $$ where $c=\cosh(t), s=\sinh(t)$ for some $t\in \R$, $\lambda\in\F$ with $|\lambda|=1$ and $A\in \U(n-1;\F)$.
\end{lemma}
\noindent Let $G\subset\U(1,n;\F)$ be a subgroup. Then ${\mathcal L}(G):=\overline{H\cdot p}\cap\partial H^n(\F)$, $p\in H^n(\F)$, is called the \textbf{limit set} of $G$. It is independend of the point $p$.
\begin{lemma}[{\cite[Lemma 4.3.4]{ChenGreenberg}}]\label{le:normalSubgroupLimitSet}
Let $N$ be a normal subgroup of $G\subset \U(1,n;\F)$. Then $G$ leaves $\mathcal{L}(N)$ invariant. Furthermore if $\mathcal{L}(N)\neq\emptyset$ and the elements of $G$ do not have a common fixed point in $\partial H^n(\F)$, then $\mathcal{L}(N)=\mathcal{L}(G)$.
\end{lemma}
\begin{theorem}[{\cite[Theorem 4.4.1]{ChenGreenberg}}]\label{th:connectedsubgroup}
Let $G$ be a connected Lie subgroup of $\U(1,n;\F)$. Then one of the following is true.
\begin{itemize}
	\item[(a)] The elements of $G$ have a common fixed point in $\overline{H^n(\F)}$.
	\item[(b)] There is a proper, totally geodesic submanifold $M$ in $H^n(\F)$ such that $\dim M>1$, ${\mathcal L}(G)=\partial M = \overline{M}\cap \partial H^n(\F)$, and $G=K\cdot U^0(M)$, where $K\subset K^0(M)$ is a connected Lie subgroup.
	\item[(c)] $\F=\C$ and $G=\SU(1,n)$.
	\item[(d)] $G=\U^0(1,n;\F)$.
\end{itemize}
\end{theorem}
\subsection{About Lie subgroups of Sp(1,\textit{n})}
\begin{corollary}\label{ko:irreducibleSp(1,n)}
Let $H\subset \Sp(1,n)$ be a connected and $\HH$-irreducible Lie subgroup. Then $H$ is conjugate to one of the following groups:
\begin{itemize}
	\item[$(i)$] $\SO^0(1,n)$, $\SO^0(1,n)\cdot \U(1)$, $\SO^0(1,n)\cdot \Sp(1)$ if $n\geq 2$,
	\item[$(ii)$] $\SU(1,n)$, $\U(1,n)$,
	\item[$(iii)$] $\Sp(1,n)$,
	\item[$(iv)$] $U^0=\left\{A\in \Sp(1,1)| A\Phirm=\Phirm A\right\}\cong\Spin^0(1,3)$ with $\Phirm=${\footnotesize$\begin{pmatrix} 0 & -1\\ 1 & 0 \end{pmatrix}$} if $n=1$.
\end{itemize}
\end{corollary}
\noindent\textit{Proof:} We apply Theorem \ref{th:connectedsubgroup} to $H$. Case $(c)$ is not relevant, since $\F=\HH$. Since every point in $\overline{H^n(\HH)}$ corresponds to a quaternionic line in $\HH^{1,n}$, every fixed point of $H$ in $\overline{H^n(\HH)}$ gives us an $H$-invariant $\HH$-subspace of $\HH^{1,n}$. Hence, we can exclude case $(a)$ in Theorem \ref{th:connectedsubgroup}. The case $(d)$ gives us $H=\U(1,n;\HH)=\Sp(1,n)$ which is $(iii)$ in the Corollary. Only case $(b)$ remains for further consideration. Here $H=K\cdot U^0(M)$ for a proper totally geodesic submanifold $M\subset H^n(\HH)$ with $\dim M>1$ and $K\subset K(M)$ connected. According to Proposition \ref{prop:totallygeodesicsubmanifolds} there are four possibilities for $M$ and furthermore by Table \ref{table:Liesubgroups} we know that $H\subset \Sp(1,m)\times\Sp(n-m)$ or $H\subset U^0\times\Sp(n-1)$. In the first case there is an invariant $\HH$-subspace if $m<n$ and we get $(i)$, $(ii)$ and $(iii)$ in the corollary. If $H\subset U^0\times\Sp(n-1)$, then $H$ can act $\HH$-irreducibly if and only if $n=1$. This gives us case $(iv)$ in the Corollary and finishes the proof. $\hfill\Box$
\begin{proposition}\label{prop:NormalSubgroupOfIrreducibleGroup}
Let $H\subset\Sp(1,n)$ be an $\HH$-irreducible subgroup. Then one of the following is true.
\begin{itemize}
	\item[$(i)$] $H$ is discrete.
	\item[$(ii)$] $H^0=\U(1)\cdot\mathbbm{1}_{n+1}$ or $H^0=\Sp(1)\cdot\mathbbm{1}_{n+1}$.
	\item[$(iii)$] $H^0$ is $\HH$-irreducible.
		\item[$(iv)$] $n=1$ and $H^0$ is one of the groups $\SO^0(1,1)$, $\SO^0(1,1)\cdot\U(1)$, $\SO^0(1,1)\cdot\Sp(1)$ or $$ S=\left\{\left. e^{ibt}\begin{pmatrix} \cosh(at) & \sinh(at)\\ \sinh(at) & \cosh(at)\end{pmatrix} \right|\ t\in\R \right\}, $$
		for some non-zero real numbers $a,b$.
\end{itemize}
\end{proposition}
\noindent\textit{Proof:} We will apply Theorem \ref{th:connectedsubgroup} to $H^0$ and discuss the cases $(a)$, $(b)$, $(c)$, and $(d)$. Assume that $H$ is not discrete. Case $(c)$ is not relevant, since $\F=\HH$. If $(d)$ holds for $H^0$, then $H^0=\Sp(1,n)$ acts $\HH$-irreducibly, so we are in $(iii)$.\\
Assume now that $(a)$ holds for $H^0$, i.e.\ $H^0$ has a common fixed point in $\overline{H^n(\HH)}$. We first discuss the case when this fixed point lies in $H^n(\HH)$. This means that all elements in $H^0$ are elliptic. Let $g\in H^0$ and $F(g)$ the set of fixed points of $g$ in $H^n(\HH)$. By Proposition \ref{prop:ellipticfixedpoints}, $F(g)$ is either a singleton or a totally geodesic submanifold. Hence, $M:=\bigcap_{g\in H^0}F(g)$ is a totally geodesic submanifold and the set of all common fixed points in $H^n(\HH)$ of $H^0$. In particular, $H^0\subset K(M)$. Since $H^0$ is a normal subgroup, $M$ is preserved by $H$, i.e.\ $H\subset I(M)$. The $\HH$-irreducibility of $H$ implies that $M$ is not a singleton. Furthermore we see from Table \ref{table:Liesubgroups} that $M$ is $H^n(\R)$, $H^n(\C)$, $H^n(\HH)$ or $n=1$ and $M=H^1(\mathbbm{I})$. Hence, the possibilities for $K(M)$ are $K(M)=\left\{\pm\mathbbm{1}_{n+1}\right\}$, $\U(1)\cdot\mathbbm{1}_{n+1}$, or $\Sp(1)\cdot\mathbbm{1}_{n+1}$. By assumption $H$ is not discrete. Since $\Sp(1)$ has no two dimensional Lie subgroup, we obtain $H^0=\U(1)\cdot\mathbbm{1}_{n+1}$ or $H^0=\Sp(1)\cdot\mathbbm{1}_{n+1}$, so we are in case $(ii)$.\\
Secondly, we consider the case, when $H^0$ has no common fixed point in $H^n(\HH)$. This means that there is a common fixed point in $\partial H^n(\HH)$. Let $F\subset\partial H^n(\HH)$ be the set of common fixed points of $H^0$ on the boundary. Notice that $F$ consists of either one or two elements, since otherwise there exist common fixed points in $H^n(\HH)$. If $F$ has exactly one element, then $H$ fixes this point, since $H^0$ is a normal subgroup of $H$. But this contradicts the $\HH$-irreducibility of $H$. If $F$ has exactly two elements, then $F$ is preserved by $H$. It follows that $H$ preserves the two dimensional $\HH$-subspace spanned by the two $H^0$-invariant lightlike lines. Since $H$ acts $\HH$-irreducibly, this can only be the case if $n=1$. By Lemma \ref{le:loxodromicconjugate} we know that $H^0$ is contained in $\SO^0(1,1)\cdot\Sp(1)$. Notice that $H^0$ is not compact, since otherwise it would be contained in a maximal compact subgroup of $\Sp(1,1)$ and would have a fixed point in $H^1(\HH)$ contradicting the assumption above. Let $H^0=L\cdot R$ be the Levi decomposition. If the Levi factor $L$ is non-trivial then it is at least three dimensional and hence equals $\Sp(1)$. In that case it follows that $R=\SO^0(1,1)$ since $H^0$ is not compact, i.e.\ $H^0=\SO^0(1,1)\cdot\Sp(1)$. If the Levi factor is trivial, then $H^0$ is contained in $\SO^0(1,1)\cdot\U(1)$. If $H^0$ has dimension two, then $H^0=\SO^0(1,1)\cdot\U(1)$. Assume now, that $H^0$ has dimension one and let $\mathfrak{h}\subset\mathfrak{so}(1,1)\oplus\mathfrak{u}(1)$ be its Lie algebra. Then $\mathfrak{h}$ is spanned by a vector $v=x+y$ with $x\in\mathfrak{so}(1,1)$ and $y\in\mathfrak{u}(1)$. We have $H^0=\exp(\R\cdot v)$. Since $H^0$ is not compact, it follows $x\neq 0$. If $y=0$, then we have $H^0=\SO^0(1,1)$. Otherwise there exist non-zero real numbers $a,b$ such that $$ v= \begin{pmatrix} 0 & a\\ a & 0\end{pmatrix} + \begin{pmatrix} ib & 0\\ 0 & ib \end{pmatrix}. $$ This gives us the group $S$. Summarizing, we are in case $(iv)$.\\
Assume now that $(b)$ holds for $H^0$, i.e.\ there is a proper totally geodesic submanifold $M$ such that $H^0=KU^0(M)$ and ${\mathcal L}(H^0)=\partial M$. Using irreducibility of $H$ we can apply Lemma \ref{le:normalSubgroupLimitSet} to $H^0$ and $H$. It follows $\partial M={\mathcal L}(H^0)={\mathcal L}(H)$. Notice that $M$ is the union of all geodesics whose endpoints lie in $\partial M$. Since the elements of $H$ map geodesics to geodesics and for the hyperbolic space the geodesics are uniquely determined by their endpoints, see \cite[Proposition 2.5.1]{ChenGreenberg}, it follows that $H$ preserves $M$. By Table \ref{table:Liesubgroups} $M$ is either $H^n(\R)$ or $H^n(\C)$. Since $H^0=KU^0(M)$, we obtain that $H^0$ acts $\HH$-irreducibly, so we are in case $(iii)$. This finishes the proof. $\hfill\Box$
\section{Main result}
\subsection{Proof of the Theorem}
\begin{lemma}\label{le:InvariantForms1}
Let $n\geq 3$ and $\alpha\in \otimes^3V^*$, where $V=\HH^{1,n}$ is considered as a real vector space. If $\alpha$ is $\SO^0(1,n)$-invariant, then $\alpha=0$.
\end{lemma}
\noindent\textit{Proof:} We have $\HH^{1,n}\cong \R^{1,n}\otimes\R^4$. Since $\SO^0(1,n)$ acts trivially on $\R^4$, we just have to check the claim for $\alpha\in \otimes^3(\R^{1,n})^*$. Since every such invariant gives rise to an invariant of $\SO(n+1,\C)$ on $\otimes^3(\C^{n+1})^*$, it is sufficient to consider these invariants. Let $v_1$, $v_2$, $v_3\in\C^{n+1}$ and $W:=\mathrm{span}_\C\left\lbrace v_1,v_2,v_3\right\rbrace$. Since $n\geq 3$, there exists some $A\in\SO(n+1,\C)$ such that $A_{|W}=-\mathrm{Id}_{W}$. It follows $\alpha(v_1,v_2,v_3)=(A^*\alpha)(v_1,v_2,v_3)=(-1)^3\alpha(v_1,v_2,v_3)$. Hence, $\alpha=0$. $\hfill\Box$
\begin{lemma}\label{le:InvariantForms2}
Let $\alpha\in \otimes^3V^*$, where $V=\HH^{1,1}$ is considered as a real vector space. If $\alpha$ is invariant under one of the Lie groups in $(iv)$ of Proposition $\ref{prop:NormalSubgroupOfIrreducibleGroup}$, then $\alpha=0$.
\end{lemma}
\noindent\textit{Proof:} We have to discuss the four Lie groups of Proposition \ref{prop:NormalSubgroupOfIrreducibleGroup} $(iv)$. The claim is clear for the groups $\SO^0(1,1)\cdot\U(1)$ and $\SO^0(1,1)\cdot\Sp(1)$, since they contain $-\mathbbm{1}_2$. Now we consider the claim for the group $\SO^0(1,1)$. As in the proof of Lemma \ref{le:InvariantForms1} it is sufficient to consider $\SO(2,\C)$ and to use that $-\mathbbm{1}_2\in \SO(2)$. So we have finally to consider the group $S$. We consider its complexification 
$$ S^{\C}=\left\{\left. e^{ib\lambda}\begin{pmatrix} \cosh(a\lambda) & \sinh(a\lambda)\\ \sinh(a\lambda) & \cosh(a\lambda)\end{pmatrix} \right|\ \lambda\in\C \right\}.$$
If we set $\lambda=i\frac{\pi}{a}$, we get an element $A=-r\mathbbm{1}_2$ with $r=e^{-\frac{\pi b}{a}}$. It follows $\alpha=A^*\alpha=-r^3\alpha$. This shows $\alpha=0$. $\hfill\Box$
\begin{remark}\label{re:InvariantForms}
The elements of $\otimes^3V^*$ that are $H$-invariant are in one-to-one correspondence with the bilinear maps $V\times V$ to $V$ that are $H$-equivariant. It follows from Lemmas $\ref{le:InvariantForms1}$ and $\ref{le:InvariantForms2}$ that the corresponding bilinear maps also vanish.
\end{remark}
\noindent Now we are able to prove Theorem \ref{th:MainTheorem}.\\[3pt]
\textit{Proof of the Theorem:} Let $\lambda: H\to \GL(T_pM)$, $h\mapsto dh_p$ be the isotropy representation. We identify $H$ with $\lambda(H)$. Since $H$ preserves the metric $g$ and the almost hypercomplex structure, we can consider $H$ as a subgroup of $\Sp(1,n)$.\\
We will first discuss the cases with $\dim M\geq 16$ and $\dim M=8$. The twelve-dimensional case is special and will be discussed afterwards.\\
We consider the universal covering $\tilde{M}=\tilde{G}/H^0$. The first step in the proof is to show that $\tilde{M}$ is a hyper-K\"ahler manifold. By Hitchin's Lemma, see \cite[Lemma 6.8]{Hitchin}, this follows by showing that the three K\"ahler forms $\omega_1$, $\omega_2$, and $\omega_3$ are closed. Since $\tilde{G}$ acts transitively, it is sufficient to show that $(d\omega_\alpha)_p=0$, $\alpha=1,2,3$.\\
First we identify the tangent space $T_p\tilde{M}$ with $\HH^{1,n}$ and consider the exterior derivatives of the three K\"ahler forms at $p$ as elements of $\mathrm{\Lambda}^3(\HH^{1,n})^*$. All three forms are invariant under the action of $H^0$. Now we apply Proposition \ref{prop:NormalSubgroupOfIrreducibleGroup} and exclude first case $(i)$. The idea is to show that it follows that $G$ is abelian, which implies that the isotropy representation is trivial contradicting our assumptions.\\
Assume that $H$ is discrete. Then $G\cong G/\left\{e\right\}\to M$ is a covering of $M$ and we can identify the Lie algebra $\g$ with $T_pM\cong\HH^{1,n}$. Notice that $H$ and its Zariski closure $H^{Zar}$ acts on $\g$ by conjugacy and that both are $\HH$-irreducible. Then the Lie bracket $\left[\cdot,\cdot\right]$ at $p$ defines an anti-symmetric bilinear on $\g\cong\HH^{1,n}$ which is $H^{Zar}$-equivariant. Since $H^{Zar}$ is an algebraic group, it has finitely many connected components, see also \cite{Milnor}. This implies that $(H^{Zar})^0$ is not compact, since otherwise $H^{Zar}$ would be compact contradicting the $\HH$-irreducibility. Proposition \ref{prop:NormalSubgroupOfIrreducibleGroup} implies that $(iii)$ or $(iv)$ holds for $(H^{Zar})^0$. Now we distinguish the cases $n\geq 3$ and $n=1$. If $n\geq 3$, then $(H^{Zar})^0$ is $\HH$-irreducible. From Corollary \ref{ko:irreducibleSp(1,n)} we see that $\SO^0(1,n)\subset (H^{Zar})^0$. By Remark \ref{re:InvariantForms}, Lemma \ref{le:InvariantForms1} implies that the Lie bracket vanishes. If $n=1$, then $(iii)$ or $(iv)$ of Proposition \ref{prop:NormalSubgroupOfIrreducibleGroup} holds for $(H^{Zar})^0$. If $(iii)$ holds, from Corollary \ref{ko:irreducibleSp(1,n)} we know that $(H^{Zar})^0$ is conjugate to one of the groups $\SU(1,1)$, $\U(1,1)$, $\Sp(1,1)$, or $U^0$. Since all four groups contain $-\mathbbm{1}_2$, the Lie bracket vanishes. If $(iv)$ holds, then Lemma \ref{le:InvariantForms2} implies together with Remark \ref{re:InvariantForms} that the Lie bracket vanishes. This gives us the contradiction and shows that $H$ is not discrete.\\
Now we can easily conclude that the K\"ahler forms are closed. We just have to consider the remaining possibilities for $H^0$ in Proposition \ref{prop:NormalSubgroupOfIrreducibleGroup}. If $(ii)$ holds, i.e. $H^0=\U(1)\cdot\mathbbm{1}_2$ or $\Sp(1)\cdot\mathbbm{1}_2$, then $-\mathbbm{1}_2\in H^0$ which implies that the K\"ahler forms are closed. If $(iii)$ or $(iv)$ holds, then this follows from Lemma \ref{le:InvariantForms1} and Lemma \ref{le:InvariantForms2}, respectively. Hence, we have shown that $\tilde{M}$ is hyper-K\"ahler.\\
Next we show that $\tilde{M}$ is a reductive homogeneous space, i.e.\ there exists an $\Ad(H^0)$-invariant vector subspace $\mathfrak{m}\subset\g$ such that $\g=\mathfrak{h}\oplus\mathfrak{m}$.\\
First we consider the case $n\geq 3$. Since $H$ is not discrete, we know from Proposition \ref{prop:NormalSubgroupOfIrreducibleGroup} and Corollary \ref{ko:irreducibleSp(1,n)} that $H^0$ is one of the following groups
$$ \U(1), \ \Sp(1), \ \SO^0(1,n), \ \SO^0(1,n)\cdot\U(1), \ \SO^0(1,n)\cdot\Sp(1) $$
$$ \SU(1,n), \ \U(1,n), \ \Sp(1,n). $$
If $H^0$ is one of the compact or semi-simple groups above, $\Ad(H^0)$ acts completely reducibly on $\g$. In particular $\tilde{M}$ is a reductive homogeneous space. So there are only the cases left where $H^0$ is $\SO^0(1,n)\cdot\U(1)$ or $\U(1,n)=\SU(1,n)\cdot\U(1)$.\\
Let $\mathfrak{s}$ be either $\mathfrak{so}(1,n)$ or $\mathfrak{su}(1,n)$. Then we have $\mathfrak{h}=\mathfrak{s}\oplus\mathfrak{u}(1)$. We consider the adjoint representation of $\mathfrak{s}$ on $\g$. Since $\mathfrak{s}$ is simple, $\mathfrak{s}$ acts completely reducibly on $\g$ and $\mathfrak{s}$ is an irreducible $\mathfrak{s}$-invariant subspace. Furthermore, there exists an $\mathfrak{s}$-invariant complement $\mathfrak{m}$ of $\mathfrak{h}=\mathfrak{s}\oplus\mathfrak{u}(1)$ which is isomorphic to $\g/\mathfrak{h}\cong T_p\tilde{M}\cong\HH^{1,n}$. Hence, the $\mathfrak{s}$-module $\g$ decomposes into $\g=\mathfrak{m}\oplus\mathfrak{s}\oplus\mathfrak{u}(1)$. Notice that $\mathfrak{m}\cong\HH^{1,n}$ decomposes into four respectively two irreducible $\mathfrak{s}$-invariant subspaces which are equivalent to $\R^{1,n}$ or $\C^{1,n}$, respectively. These three submodules $\mathfrak{s}$, $\mathfrak{u}(1)$, $\R^{1,n}$ or $\C^{1,n}$ are pairwise inequivalent. Since $\mathfrak{s}$ and $\mathfrak{u}(1)$ commute, $\mathfrak{u}(1)$ preserves the isotypical $\mathfrak{s}$-submodules. It follows that the isotypical submodule $\mathfrak{m}$ is $\mathfrak{u}(1)$-invariant and thus also $\mathfrak{h}$-invariant. Hence, $\mathfrak{m}$ is invariant under $\Ad(H^0)$.\\
Now we investigate the case $n=1$. By Proposition \ref{prop:NormalSubgroupOfIrreducibleGroup} and Corollary \ref{ko:irreducibleSp(1,n)} we know that $H^0$ is one of the following groups
$$ \U(1), \ \Sp(1), \ \SU(1,1), \ \U(1,1), \ \Sp(1,1), \ U^0\cong \Spin^0(1,3), $$
$$ \SO^0(1,1)\cdot\Sp(1), \ \SO^0(1,1)\cdot\U(1), \ \SO^0(1,1), \ S. $$
If $H^0$ is one of the compact or semi-simple groups it is clear that $\tilde{M}$ is a reductive homogeneous space. If $H^0$ is $\U(1,1)$ then we can apply the arguments from above. So we just have to consider the Lie groups $\SO^0(1,1)\cdot\Sp(1)$, $\SO^0(1,1)\cdot\U(1)$, $\SO^0(1,1)$, and $S$. Recall that $\g\cong \mathfrak{h}\oplus\g/\mathfrak{h}$ and $\g/\mathfrak{h}\cong T_p\tilde{M}\cong \HH^{1,1}$.\\[3pt]
If $H^0=\SO^0(1,1)\cdot\Sp(1)$ then $\g$ decomposes in $\Sp(1)$-submodules, namely $\g=\mathfrak{so}(1,1)\oplus\mathfrak{sp}(1)\oplus V\oplus V\cong \R\oplus\R^3\oplus V\oplus V$ with $\dim_\R V=4$. Let $W:=V\oplus V\cong\R^8$. Hence, we have three isotypical submodules which are $\Sp(1)$-invariant. Since the elements of $\SO^0(1,1)$ and $\Sp(1)$ are commuting, $\SO^0(1,1)$ preserves this decomposition into isotypical submodules. In particular $\SO^0(1,1)$ preserves the complement $W$ of $\mathfrak{h}$.\\[3pt]
If $H^0=\SO^0(1,1)\cdot\U(1)$, then $\g$ decomposes into $\U(1)$-invariant submodules $\g=\mathfrak{so}(1,1)\oplus\mathfrak{u}(1)\oplus(V\oplus V\oplus V\oplus V)$ with $V\cong \C$. So there are two isotypical submodules. As before it follows that $\mathfrak{m}=\oplus^4 V$ is an $H^0$-invariant complement of $\mathfrak{h}\subset\g$.\\[3pt]
Assume $H^0=\SO^0(1,1)$ and consider the adjoint action of $\mathfrak{h}=\mathfrak{so}(1,1)$ on $\g\cong\R^9$. Let $A\in\mathfrak{so}(1,1)\setminus\left\{0\right\}$. Then $\g$ decomposes into $\g=\ker A\oplus \mathrm{im } \ A$, where $\mathrm{im } \ A=V_+\oplus V_-$ is a sum of two $4$-dimensional eigenspaces with opposite real  eigenvalues
and provides an $\Ad(H^0)$-invariant complement to $\mathfrak{h}\subset\g$.\\[3pt]
If $H^0=S$ then $\g$ decomposes as before into $\g=\ker A\oplus \textrm{im }A$ where $0\neq A\in \mathfrak{h}$ and $\textrm{im }A$ is $8$-dimensional (a sum 
of $2$ complex eigenspaces of equal dimension). As before this gives us an $\Ad(H^0)$-invariant complement to $\mathfrak{h}\subset \g$.\\[3pt]
Summarizing we have shown that $\tilde{M}$ is indeed a reductive homogeneous space. Next we show that $\g=\mathfrak{h}\oplus\mathfrak{m}$ is a symmetric Lie algebra. It is sufficient to show that $\left[\mathfrak{m},\mathfrak{m} \right]\subset\mathfrak{h}$. We consider the restriction of the Lie bracket $\left[\cdot,\cdot\right]$ to $\mathfrak{m}\times\mathfrak{m}$ and denote its projection to $\mathfrak{m}$ by $\beta$. The antisymmetric bilinear map $\beta$ is $\Ad(H^0)$-equivariant. Since $\mathfrak{m}\cong \HH^{1,n}$, we can consider such a map as an $H^0$-invariant element of $\otimes^3(\HH^{1,n})^*$. We already discussed that such maps vanish, so we have $\beta=0$. This proves $\left[\mathfrak{m},\mathfrak{m} \right]\subset\mathfrak{h}$.\\
Next we show that $\tilde{M}$ is isometric to $\Mink_{n+1}(\HH)$. We consider the Lie algebra of the transvection group $\hat{\g}=\left[\mathfrak{m},\mathfrak{m}\right]\oplus\mathfrak{m}\subset\g$. It is know that the transvection group of a hyper-K\"ahler symmetric space is nilpotent, see \cite[Corollary 2.4]{KathOlbrich}. It follows that the action of $\hat{\g}_+:=\left[\mathfrak{m},\mathfrak{m}\right]\subset\mathfrak{h}$ on $\mathfrak{m}$ is nilpotent. From the above list of the possible Lie algebras we see that $\mathfrak{h}$ is reductive, i.e.\ $\mathfrak{h}=\mathfrak{s}\oplus\mathfrak{a}$ where $\mathfrak{s}$ is semi-simple and $\mathfrak{a}$ is abelian. Notice that $\mathfrak{s}$ is also allowed to be trivial. Since $\hat{\g}_+$ is nilpotent and furthermore an ideal in $\mathfrak{h}$, it follows that $\hat{\g}_+\subset \mathfrak{a}$. Since $\mathfrak{a}$ acts completely reducibly on $\mathfrak{m}$, the same holds for $\hat{\g}_+$. Hence, $\hat{\g}_+$ acts trivially on $\mathfrak{m}$. Since the action of $\hat{\g}_+$ is faithful, it follows that $\hat{\g}_+=\left\lbrace 0\right\rbrace$. This proves that $\tilde{M}$ is isometric to $\Mink_{n+1}(\HH)$.\\[3pt]
Finally, we discuss the twelve-dimensional case. In this situation $H$ is discrete or $H^0$ is one of the following groups
$$ \U(1), \ \Sp(1), \ \SO^0(1,2), \ \SO^0(1,2)\cdot\U(1), \ \SO^0(1,2)\cdot\Sp(1),$$
$$ \ \SU(1,2), \ \U(1,2), \ \Sp(1,2).$$
Notice that there exist non-trivial 3-forms that are invariant under $\SO^0(1,2)$ or $\SU(1,2)$. In particular we can not exclude that $H$ is discrete since $(H^{Zar})^0$ could be $\SO^0(1,2)$. But if $H^0$ is one of the groups $\U(1)$, $\Sp(1)$, $\SO^0(1,2)\cdot\U(1)$, $\SO^0(1,2)\cdot\Sp(1)$, $\U(1,2)$, or $\Sp(1,2)$, then we can apply all arguments from above and get $\tilde{M}\cong \Mink_{n+1}(\HH)$. Thus if $\tilde{M}$ is not isometric to $\Mink_{n+1}(\HH)$, then we get from the above list that $H^0$ is either $\left\{e\right\}$, $\SO^0(1,2)$, or $\SU(1,2)$. This finishes the proof. $\hfill\Box$

\subsection{A class of non-symmetric examples in dimension 12}
By Theorem \ref{th:MainTheorem} we know that a non-flat manifold appears only if $\dim M=12$ and one of the following is true
\begin{itemize}
 \item $\mathfrak{h}=\mathfrak{so}(1,2)$,
 \item $\mathfrak{h}=\mathfrak{su}(1,2)$,
 \item $\mathfrak{h}=\left\lbrace 0\right\rbrace$, but $(H^{Zar})^0=\SO^0(1,2)$ or $\SU(1,2)$.
\end{itemize}
This is due to the fact that there exist non-trivial $3$-forms on $\HH^{1,2}$ which are invariant under $\SO^0(1,2)$. Therefore we can not conlude that the K\"ahler forms are closed if $\SO^0(1,2)\subset H^0$ or $\SO^0(1,2)\subset H^{Zar}$.\\
In the following we will investigate the case with $\mathfrak{h}=\mathfrak{so}(1,2)$ and give some non-symmetric examples. We consider the following Lie algebra
$$ \mathfrak{m}=\ell\cdot\mathfrak{so}(1,2)\oplus \R^{1,2}\otimes\R^{4-\ell}$$
for $\ell\in\left\lbrace 1,2,3,4\right\rbrace$ and where the subalgebra $\R^{1,2}\otimes\R^{4-\ell}$ is abelian. We define the representation $\rho:\mathfrak{so}(1,2)\to\mathrm{der }(\mathfrak{m})$ by the adjoint representation on $\ell\cdot\mathfrak{so}(1,2)$, by the standard representation on $\R^{1,2}$ and by the trivial representation on $\R^{4-\ell}$.\\
Now we set $\mathfrak{h}=\mathfrak{so}(1,2)$ and $\g=\mathfrak{h}\ltimes_{\rho}\mathfrak{m}\cong (\ell+1)\cdot\mathfrak{so}(1,2)\ltimes_{\rho'}\R^{1,2}\otimes\R^{4-\ell}$ where $\mathfrak{h}$ corresponds to $\lbrace (X,X,\ldots,X)\in (\ell+1)\cdot\mathfrak{so}(1,2)\ | \ X\in\mathfrak{so}(1,2) \rbrace$ and $\rho'(X_0,\ldots,X_\ell)x\otimes v = X_0x\otimes v$ for all $x\in\R^{1,2}$, $v\in\R^{4-\ell}$, $X_0,\ldots,X_\ell\in\mathfrak{so}(1,2)$.\\
The isotropy representation is equivalent to $\R^{1,2}\otimes\R^4\cong\HH^{1,2}$, hence, admits an $\mathfrak{h}$-invariant hyper-Hermitian structure of index 4.\\[3pt]
For a general classification of the homogeneous spaces with $\mathfrak{h}=\mathfrak{so}(1,2)$ one needs to classify all Lie algebra structures on the vector space $\g=\mathfrak{so}(1,2)\oplus\R^{1,2}\otimes\R^4$ such that the Lie bracket restricts to the Lie bracket of $\mathfrak{so}(1,2)$ and to the canonical representation of $\mathfrak{so}(1,2)$ on $\R^{1,2}\otimes\R^4$. For this one has to describe all $\mathfrak{so}(1,2)$-invariant tensors in $\mathrm{\Lambda}^2(\R^{1,2}\otimes\R^4)^*\otimes\g\cong \mathrm{\Lambda}^2(\R^{1,2}\otimes\R^4)\otimes(\R^{1,2}\otimes\R^5)$ which satisfy the Jacobi identity. Since
$$ \left[\mathrm{\Lambda}^2(\R^{1,2}\otimes\R^4)\otimes (\R^{1,2}\otimes\R^5) \right]^{\mathfrak{so}(1,2)}\cong \mathrm{\Lambda}^3\R^{1,2}\otimes S^2\R^4\otimes\R^5 $$
the Lie brackets on $\mathfrak{m}=\R^{1,2}\otimes\R^4$ are of the form
$$ \left[ x\otimes v, y\otimes w\right] = (x\times y)\otimes \beta(v,w), $$
where $\beta\in S^2(\R^4)^*\otimes\R^5$, so $\beta=\sum_{i=0}^4\beta_i\otimes b_i$, where $(b_i)$ is a basis of $\R^5$.
The Jacobi identity for three vectors in $\mathfrak{g}$ then holds if at least one of the three vectors is in $\mathfrak{h}$ and the remaining equations form a
system of quadratic equations for $\beta$. The above examples correspond to 
solutions of the form $\beta_0=0$ and $\beta_i=\lambda_i(b_i^*)^2$, $i=1,2,3,4$.\\
Now we consider the intrinsic torsion of the indefinite almost hyper-Hermitian structure on $M$. From the fact that 
the sum $\nabla +S$ of the Levi-Civita connection $\nabla$ and the tensor field $S=-\frac{1}{4}\sum_{\alpha =1}^3J_\alpha\nabla J_\alpha$ is a connection compatible with the metric and the 
hypercomplex structure, it follows that the intrinsic torsion is completely determined by the three covariant derivatives $\nabla 
J_\alpha$. More precisely, it is given by the image of $S$ (evaluated at the canonical base point) under the isomorphism $V^*\otimes \mathfrak{sp}(1,2)^\perp\cong \frac{\mathrm{\Lambda}^2V^*\otimes V}{\mathrm{alt}(V^*\otimes\, \mathfrak{sp}(1,2))}$ induced by the alternation map 
$\mathrm{alt} : V^*\otimes \mathfrak{so}(V) \stackrel{\sim}{\longrightarrow} \mathrm{\Lambda}^2V^*\otimes V$, where $V=\R^{12}$ and 
$\perp$ stands for the orthogonal complement in $\mathfrak{so}(V)$ with respect to the Killing form.

We have the following formula for each almost complex structure
$$ 4g((\nabla_XJ_\alpha)Y,Z) = 6d\omega_\alpha(X,J_\alpha Y,J_\alpha Z)-6d\omega_\alpha(X,Y,Z)+g(N_{J_\alpha}(Y,Z),J_\alpha X), $$
where $N_{J_\alpha}$ denotes the Nijenhuis tensor, see \cite{KobayashiNomizu}. Therefore, $\nabla J_\alpha$ is determined by $d\omega_\alpha$ and $N_{J_\alpha}$. The 3-form $d\omega_\alpha$ is $\SO^0(1,2)$-invariant. Since $$ (\mathrm{\Lambda}^3(\R^{1,2}\otimes\R^4))^{\SO^0(1,2)}\cong \mathrm{\Lambda}^3 \R^{1,2}\otimes S^3\R^4, $$
it has the form $d\omega_\alpha(x_1\otimes q_1,x_2\otimes q_2,x_3\otimes q_3) = s_\alpha \cdot\det(x_1,x_2,x_3)\sigma_\alpha (q_1,q_2,q_3)$ with $x_i\otimes q_i\in\R^{1,2}\otimes\R^4$, $s_\alpha\in\R$, and $\sigma_\alpha\in S^3(\R^4)^*$, where $\R^{1,2}\otimes\R^4$ is identified with the tangent space of $M$ at the canonical base point. Analogously, the Nijenhuis tensor is given by
$$ N_{J_\alpha}(x_1\otimes q_1,x_2\otimes q_2) = t_\alpha\cdot K(x_1,x_2)\tau_\alpha(q_1,q_2)$$
where $x_1,x_2\in\R^{1,2}$, $q_1,q_2\in\R^4$, $t_\alpha \in\R$, $\tau_\alpha\in S^2(\R^4)^*\otimes\R^4$ and $K$ is the cross product on $\R^{1,2}$, since $(\mathrm{\Lambda}^2(\R^{1,2}\otimes\R^4)^*\otimes\R^{1,2}\otimes\R^4)^{\SO^0(1,2)}\cong \R\cdot K\otimes S^2(\R^4)^*\otimes\R^4$.\\

Finally we give an example for the case $\mathfrak{h}=0$ but $(H^{Zar})^0=\SO^0(1,2)$. As before, let $\mathfrak{m}=\ell\cdot\mathfrak{so}(1,2)\oplus\R^{1,2}\otimes\R^{4-\ell}$ with $\ell\in\left\{1,2,3,4\right\}$. Let $H$ be the image of $\SL(2,\Z)$ under the double cover $\SL(2,\R)\to\SO^0(1,2)$. Then $H$ is a discrete, Zariski dense and $\HH$-irreducible subgroup. If we set $G=H\ltimes \left( \SO^0(1,2)^\ell \times \R^{3(4-\ell)}\right)$, we get the desired homogeneous space $M=G/H$.\\

\end{document}